\documentclass{gtart}
\def\ifplaintex{\expandafter\ifx\csname documentclass\endcsname\relax}

\def\gtp{{\mathsurround=0pt\it $\cal G\mskip-2mu$eometry \&\ 
$\cal T\!\!$opology $\cal P\!$ublications}}  % GT publications

\def\recd{{\small Received:\qua\receiveddate\ifx\reviseddate\relax
\else\qquad Revised:\qua\reviseddate\fi\par}} 

\def\lognumber#1{\def\thelognumber{#1}}
\def\volumenumber#1{\def\thevolumenumber{#1}}
\def\volumeyear#1{\def\thevolumeyear{#1}}
\def\papernumber#1{\def\thepapernumber{#1}}
\def\pagenumbers#1#2{\def\startpage{#1}\def\finishpage{#2}}
\def\published#1{\def\publishdate{#1}}

\def\received#1{\def\receiveddate{#1}}
\def\revised#1{\def\reviseddate{#1}}
\def\accepted#1{\def\accepteddate{#1}}

\def\asciiaddress#1{\def\theasciiaddress{#1}}

\long\def\asciiabstract#1{\long\def\theasciiabstract{#1}}

\let\\\par\let\thelognumber\relax\let\thevolumenumber\relax
\let\thepapernumber\relax\let\thevolumeyear\relax\let\startpage\relax
\let\finishpage\relax\let\publishdate\relax\let\receiveddate\relax
\let\reviseddate\relax\let\accepteddate\relax\let\theasciititle\relax
\let\theasciiauthors\relax\let\theasciiaddress\relax
\let\theasciiabstract\relax

\let\theasciiemail\relax

\ifplaintex
\font\logobig=cmssbx10 scaled 3836
\font\logomed=cmssbx10 scaled 2557
\else
\font\logobig=cmssbx10 scaled 4200
\font\logomed=cmssbx10 scaled 2800
\fi

\long\def\makeagttitle{   %%% start of definition of \makeagttitle
\count0=\startpage
\agt\hfill      %   Journal title (top left) 
\hbox to 45truept{\vbox to 0pt{\vglue -13truept{\logomed A\kern -.37em{\logobig 
T}\kern -.38em G}\vss}\hss}
\break
{\small Volume \thevolumenumber\ (\thevolumeyear)
\startpage--\finishpage\nl
Published: \publishdate}

\vglue .25truein

{\parskip=0pt\leftskip 0pt plus
1fil\def\\{\par\smallskip}{\Large\bf\thetitle}\par\medskip} \vglue
0.05truein

{\parskip=0pt\leftskip 0pt plus 1fil\def\\{\par}{\sc\theauthors}
\par\medskip}%
 
\vglue 0.03truein 

{\small\leftskip 25truept\rightskip 25truept{\bf Abstract}\stdspace\theabstract

{\bf AMS Classification}\stdspace\theprimaryclass
\ifx\thesecondaryclass\relax\else; \thesecondaryclass\fi\par
{\bf Keywords}\stdspace \thekeywords\par}\vglue 7truept

}   %%%% end of definition of \makeagttitle

\ifplaintex
\font\phead=cmsl9 scaled 950
\font\pnum=cmbx10 scaled 913
\font\pfoot=cmsl9 scaled 950
\headline{\vbox to 0pt{\vskip -4.5mm\line{\small\phead\ifnum
\count0=\startpage ISSN 1472-2739 (on-line) 1472-2747 (printed)
\hfill {\pnum\folio}\else\ifodd\count0\def\\{ }% 
\ifx\theshorttitle\relax\thetitle\else\theshorttitle\fi\hfill{\pnum\folio}
\else\def\\{ and }{\pnum\folio}\hfill\ifx\theshortauthors\relax\theauthors
\else\theshortauthors\fi\fi\fi}\vss}}
\footline{\vbox to 0pt{\vglue 0mm\line{\small\pfoot\ifnum\count0=\startpage
\copyright\ \gtp\hfill\else
\agt, Volume \thevolumenumber\ (\thevolumeyear)\hfill\fi}\vss}}
\else
\font=cmsl9 scaled 1050
\font\lnum=cmbx10 
\font=cmsl9 scaled 1050
\makeatletter
\def\@oddhead{{\small\lhead\ifnum\count0=\startpage ISSN 1472-2739 
(on-line) 1472-2747 (printed)\hfill {\lnum\number\count0}\else\ifodd\count0
\def\\{ }\ifx\theshorttitle\relax \thetitle \else\theshorttitle\fi\hfill
{\lnum\number\count0}\else\def\\{ and }{\lnum\number\count0}
\hfill\ifx\theshortauthors\relax 
\theauthors\else\theshortauthors\fi\fi\fi}}\def\@evenhead{\@oddhead}
\def\@oddfoot{\small\lfoot\ifnum\count0=\startpage\copyright\ \gtp\hfill\else
\agt, Volume \thevolumenumber\ (\thevolumeyear)\hfill\fi}
\def\@evenfoot{\@oddfoot}
\makeatother
\fi
\let\maketitlepage\makeagttitle
\let\makeshorttitle\maketitlepage
\let\maketitle\maketitlepage

\newwrite\gtoutfile
\long\gdef\makeheadfile{  %%% start of definition of \makeheadfile
{\def\\{, }\def\s{ }
\immediate\openout\gtoutfile head.xxx
\immediate\write\gtoutfile{To: math@arxiv.org}
\immediate\write\gtoutfile{Subject: put OR rep NNNNN:ppppp}
\immediate\write\gtoutfile{--text follows this line--}
\immediate\write\gtoutfile{Proxy-for: \ifx\theasciiauthors\relax
\theauthors\else\theasciiauthors\fi\s<\ifx\theasciiemail\relax\theemail\else\theasciiemail\fi>}
\immediate\write\gtoutfile{\noexpand\\}
\immediate\write\gtoutfile{Authors: \ifx\theasciiauthors\relax
\theauthors\else\theasciiauthors\fi}
{\def\\{ }\immediate\write\gtoutfile{Title: \ifx\theasciititle\relax
\thetitle\else\theasciititle\fi}}
\immediate\write\gtoutfile{Subj-class: GT or SG, GR etc}
\immediate\write\gtoutfile{MSC-class: \theprimaryclass\ifx\thesecondaryclass\relax\else, \thesecondaryclass\fi}
\immediate\write\gtoutfile{Journal-ref: Algebr. Geom. Topol. \thevolumenumber\s
(\thevolumeyear) \startpage-\finishpage}
\immediate\write\gtoutfile{Comments: Published by Algebraic and
Geometric Topology at}
\immediate\write\gtoutfile{\s\s\s  http://www.maths.warwick.ac.uk/agt/AGTVol\thevolumenumber/agt-\thevolumenumber-\thepapernumber.abs.html}
\immediate\write\gtoutfile{\noexpand\\}
\immediate\write\gtoutfile{}
\ifx\theasciiabstract\relax
\immediate\write\gtoutfile{\theabstract}\else
\immediate\write\gtoutfile{\theasciiabstract}\fi
\immediate\write\gtoutfile{}
\immediate\write\gtoutfile{\noexpand\\}
\immediate\write\gtoutfile{}
\immediate\closeout\gtoutfile}}  %%% end of definition of \makeheadfile

\def\maketitlepage{\makeagttitle\makeheadfile}
\let\makeshorttitle\maketitlepage
\let\maketitle\maketitlepage

\def\ifplaintex{\expandafter\ifx\csname documentclass\endcsname\relax}

\def\gtp{{\mathsurround=0pt\it $\cal G\mskip-2mu$eometry \&\ 
$\cal T\!\!$opology $\cal P\!$ublications}}  % GT publications

\def\recd{{\small Received:\qua\receiveddate\ifx\reviseddate\relax
\else\qquad Revised:\qua\reviseddate\fi\par}} 

\def\lognumber#1{\def\thelognumber{#1}}
\def\volumenumber#1{\def\thevolumenumber{#1}}
\def\volumeyear#1{\def\thevolumeyear{#1}}
\def\papernumber#1{\def\thepapernumber{#1}}
\def\pagenumbers#1#2{\def\startpage{#1}\def\finishpage{#2}}
\def\published#1{\def\publishdate{#1}}

\def\received#1{\def\receiveddate{#1}}
\def\revised#1{\def\reviseddate{#1}}
\def\accepted#1{\def\accepteddate{#1}}

\def\asciiaddress#1{\def\theasciiaddress{#1}}

\long\def\asciiabstract#1{\long\def\theasciiabstract{#1}}

\let\\\par\let\thelognumber\relax\let\thevolumenumber\relax
\let\thepapernumber\relax\let\thevolumeyear\relax\let\startpage\relax
\let\finishpage\relax\let\publishdate\relax\let\receiveddate\relax
\let\reviseddate\relax\let\accepteddate\relax\let\theasciititle\relax
\let\theasciiauthors\relax\let\theasciiaddress\relax
\let\theasciiabstract\relax

\let\theasciiemail\relax

\ifplaintex
\font\logobig=cmssbx10 scaled 3836
\font\logomed=cmssbx10 scaled 2557
\else
\font\logobig=cmssbx10 scaled 4200
\font\logomed=cmssbx10 scaled 2800
\fi

\long\def\makeagttitle{   %%% start of definition of \makeagttitle
\count0=\startpage
\agt\hfill      %   Journal title (top left) 
\hbox to 45truept{\vbox to 0pt{\vglue -13truept{\logomed A\kern -.37em{\logobig 
T}\kern -.38em G}\vss}\hss}
\break
{\small Volume \thevolumenumber\ (\thevolumeyear)
\startpage--\finishpage\nl
Published: \publishdate}

\vglue .25truein

{\parskip=0pt\leftskip 0pt plus
1fil\def\\{\par\smallskip}{\Large\bf\thetitle}\par\medskip} \vglue
0.05truein

{\parskip=0pt\leftskip 0pt plus 1fil\def\\{\par}{\sc\theauthors}
\par\medskip}%
 
\vglue 0.03truein 

{\small\leftskip 25truept\rightskip 25truept{\bf Abstract}\stdspace\theabstract

{\bf AMS Classification}\stdspace\theprimaryclass
\ifx\thesecondaryclass\relax\else; \thesecondaryclass\fi\par
{\bf Keywords}\stdspace \thekeywords\par}\vglue 7truept

}   %%%% end of definition of \makeagttitle

\ifplaintex
\font\phead=cmsl9 scaled 950
\font\pnum=cmbx10 scaled 913
\font\pfoot=cmsl9 scaled 950
\headline{\vbox to 0pt{\vskip -4.5mm\line{\small\phead\ifnum
\count0=\startpage ISSN 1472-2739 (on-line) 1472-2747 (printed)
\hfill {\pnum\folio}\else\ifodd\count0\def\\{ }% 
\ifx\theshorttitle\relax\thetitle\else\theshorttitle\fi\hfill{\pnum\folio}
\else\def\\{ and }{\pnum\folio}\hfill\ifx\theshortauthors\relax\theauthors
\else\theshortauthors\fi\fi\fi}\vss}}
\footline{\vbox to 0pt{\vglue 0mm\line{\small\pfoot\ifnum\count0=\startpage
\copyright\ \gtp\hfill\else
\agt, Volume \thevolumenumber\ (\thevolumeyear)\hfill\fi}\vss}}
\else
\font=cmsl9 scaled 1050
\font\lnum=cmbx10 
\font=cmsl9 scaled 1050
\makeatletter
\def\@oddhead{{\small\lhead\ifnum\count0=\startpage ISSN 1472-2739 
(on-line) 1472-2747 (printed)\hfill {\lnum\number\count0}\else\ifodd\count0
\def\\{ }\ifx\theshorttitle\relax \thetitle \else\theshorttitle\fi\hfill
{\lnum\number\count0}\else\def\\{ and }{\lnum\number\count0}
\hfill\ifx\theshortauthors\relax 
\theauthors\else\theshortauthors\fi\fi\fi}}\def\@evenhead{\@oddhead}
\def\@oddfoot{\small\lfoot\ifnum\count0=\startpage\copyright\ \gtp\hfill\else
\agt, Volume \thevolumenumber\ (\thevolumeyear)\hfill\fi}
\def\@evenfoot{\@oddfoot}
\makeatother
\fi
\let\maketitlepage\makeagttitle
\let\makeshorttitle\maketitlepage
\let\maketitle\maketitlepage

\newwrite\gtoutfile
\long\gdef\makeheadfile{  %%% start of definition of \makeheadfile
{\def\\{, }\def\s{ }
\immediate\openout\gtoutfile head.xxx
\immediate\write\gtoutfile{To: math@arxiv.org}
\immediate\write\gtoutfile{Subject: put OR rep NNNNN:ppppp}
\immediate\write\gtoutfile{--text follows this line--}
\immediate\write\gtoutfile{Proxy-for: \ifx\theasciiauthors\relax
\theauthors\else\theasciiauthors\fi\s<\ifx\theasciiemail\relax\theemail\else\theasciiemail\fi>}
\immediate\write\gtoutfile{\noexpand\\}
\immediate\write\gtoutfile{Authors: \ifx\theasciiauthors\relax
\theauthors\else\theasciiauthors\fi}
{\def\\{ }\immediate\write\gtoutfile{Title: \ifx\theasciititle\relax
\thetitle\else\theasciititle\fi}}
\immediate\write\gtoutfile{Subj-class: GT or SG, GR etc}
\immediate\write\gtoutfile{MSC-class: \theprimaryclass\ifx\thesecondaryclass\relax\else, \thesecondaryclass\fi}
\immediate\write\gtoutfile{Journal-ref: Algebr. Geom. Topol. \thevolumenumber\s
(\thevolumeyear) \startpage-\finishpage}
\immediate\write\gtoutfile{Comments: Published by Algebraic and
Geometric Topology at}
\immediate\write\gtoutfile{\s\s\s  http://www.maths.warwick.ac.uk/agt/AGTVol\thevolumenumber/agt-\thevolumenumber-\thepapernumber.abs.html}
\immediate\write\gtoutfile{\noexpand\\}
\immediate\write\gtoutfile{}
\ifx\theasciiabstract\relax
\immediate\write\gtoutfile{\theabstract}\else
\immediate\write\gtoutfile{\theasciiabstract}\fi
\immediate\write\gtoutfile{}
\immediate\write\gtoutfile{\noexpand\\}
\immediate\write\gtoutfile{}
\immediate\closeout\gtoutfile}}  %%% end of definition of \makeheadfile

\def\maketitlepage{\makeagttitle\makeheadfile}
\let\makeshorttitle\maketitlepage
\let\maketitle\maketitlepage

\def\ifplaintex{\expandafter\ifx\csname documentclass\endcsname\relax}

\def\gtp{{\mathsurround=0pt\it $\cal G\mskip-2mu$eometry \&\ 
$\cal T\!\!$opology $\cal P\!$ublications}}  % GT publications

\def\recd{{\small Received:\qua\receiveddate\ifx\reviseddate\relax
\else\qquad Revised:\qua\reviseddate\fi\par}} 

\def\lognumber#1{\def\thelognumber{#1}}
\def\volumenumber#1{\def\thevolumenumber{#1}}
\def\volumeyear#1{\def\thevolumeyear{#1}}
\def\papernumber#1{\def\thepapernumber{#1}}
\def\pagenumbers#1#2{\def\startpage{#1}\def\finishpage{#2}}
\def\published#1{\def\publishdate{#1}}

\def\received#1{\def\receiveddate{#1}}
\def\revised#1{\def\reviseddate{#1}}
\def\accepted#1{\def\accepteddate{#1}}

\def\asciiaddress#1{\def\theasciiaddress{#1}}

\long\def\asciiabstract#1{\long\def\theasciiabstract{#1}}

\let\\\par\let\thelognumber\relax\let\thevolumenumber\relax
\let\thepapernumber\relax\let\thevolumeyear\relax\let\startpage\relax
\let\finishpage\relax\let\publishdate\relax\let\receiveddate\relax
\let\reviseddate\relax\let\accepteddate\relax\let\theasciititle\relax
\let\theasciiauthors\relax\let\theasciiaddress\relax
\let\theasciiabstract\relax

\let\theasciiemail\relax

\ifplaintex
\font\logobig=cmssbx10 scaled 3836
\font\logomed=cmssbx10 scaled 2557
\else
\font\logobig=cmssbx10 scaled 4200
\font\logomed=cmssbx10 scaled 2800
\fi

\long\def\makeagttitle{   %%% start of definition of \makeagttitle
\count0=\startpage
\agt\hfill      %   Journal title (top left) 
\hbox to 45truept{\vbox to 0pt{\vglue -13truept{\logomed A\kern -.37em{\logobig 
T}\kern -.38em G}\vss}\hss}
\break
{\small Volume \thevolumenumber\ (\thevolumeyear)
\startpage--\finishpage\nl
Published: \publishdate}

\vglue .25truein

{\parskip=0pt\leftskip 0pt plus
1fil\def\\{\par\smallskip}{\Large\bf\thetitle}\par\medskip} \vglue
0.05truein

{\parskip=0pt\leftskip 0pt plus 1fil\def\\{\par}{\sc\theauthors}
\par\medskip}%
 
\vglue 0.03truein 

{\small\leftskip 25truept\rightskip 25truept{\bf Abstract}\stdspace\theabstract

{\bf AMS Classification}\stdspace\theprimaryclass
\ifx\thesecondaryclass\relax\else; \thesecondaryclass\fi\par
{\bf Keywords}\stdspace \thekeywords\par}\vglue 7truept

}   %%%% end of definition of \makeagttitle

\ifplaintex
\font\phead=cmsl9 scaled 950
\font\pnum=cmbx10 scaled 913
\font\pfoot=cmsl9 scaled 950
\headline{\vbox to 0pt{\vskip -4.5mm\line{\small\phead\ifnum
\count0=\startpage ISSN 1472-2739 (on-line) 1472-2747 (printed)
\hfill {\pnum\folio}\else\ifodd\count0\def\\{ }% 
\ifx\theshorttitle\relax\thetitle\else\theshorttitle\fi\hfill{\pnum\folio}
\else\def\\{ and }{\pnum\folio}\hfill\ifx\theshortauthors\relax\theauthors
\else\theshortauthors\fi\fi\fi}\vss}}
\footline{\vbox to 0pt{\vglue 0mm\line{\small\pfoot\ifnum\count0=\startpage
\copyright\ \gtp\hfill\else
\agt, Volume \thevolumenumber\ (\thevolumeyear)\hfill\fi}\vss}}
\else
\font=cmsl9 scaled 1050
\font\lnum=cmbx10 
\font=cmsl9 scaled 1050
\makeatletter
\def\@oddhead{{\small\lhead\ifnum\count0=\startpage ISSN 1472-2739 
(on-line) 1472-2747 (printed)\hfill {\lnum\number\count0}\else\ifodd\count0
\def\\{ }\ifx\theshorttitle\relax \thetitle \else\theshorttitle\fi\hfill
{\lnum\number\count0}\else\def\\{ and }{\lnum\number\count0}
\hfill\ifx\theshortauthors\relax 
\theauthors\else\theshortauthors\fi\fi\fi}}\def\@evenhead{\@oddhead}
\def\@oddfoot{\small\lfoot\ifnum\count0=\startpage\copyright\ \gtp\hfill\else
\agt, Volume \thevolumenumber\ (\thevolumeyear)\hfill\fi}
\def\@evenfoot{\@oddfoot}
\makeatother
\fi
\let\maketitlepage\makeagttitle
\let\makeshorttitle\maketitlepage
\let\maketitle\maketitlepage

\newwrite\gtoutfile
\long\gdef\makeheadfile{  %%% start of definition of \makeheadfile
{\def\\{, }\def\s{ }
\immediate\openout\gtoutfile head.xxx
\immediate\write\gtoutfile{To: math@arxiv.org}
\immediate\write\gtoutfile{Subject: put OR rep NNNNN:ppppp}
\immediate\write\gtoutfile{--text follows this line--}
\immediate\write\gtoutfile{Proxy-for: \ifx\theasciiauthors\relax
\theauthors\else\theasciiauthors\fi\s<\ifx\theasciiemail\relax\theemail\else\theasciiemail\fi>}
\immediate\write\gtoutfile{\noexpand\\}
\immediate\write\gtoutfile{Authors: \ifx\theasciiauthors\relax
\theauthors\else\theasciiauthors\fi}
{\def\\{ }\immediate\write\gtoutfile{Title: \ifx\theasciititle\relax
\thetitle\else\theasciititle\fi}}
\immediate\write\gtoutfile{Subj-class: GT or SG, GR etc}
\immediate\write\gtoutfile{MSC-class: \theprimaryclass\ifx\thesecondaryclass\relax\else, \thesecondaryclass\fi}
\immediate\write\gtoutfile{Journal-ref: Algebr. Geom. Topol. \thevolumenumber\s
(\thevolumeyear) \startpage-\finishpage}
\immediate\write\gtoutfile{Comments: Published by Algebraic and
Geometric Topology at}
\immediate\write\gtoutfile{\s\s\s  http://www.maths.warwick.ac.uk/agt/AGTVol\thevolumenumber/agt-\thevolumenumber-\thepapernumber.abs.html}
\immediate\write\gtoutfile{\noexpand\\}
\immediate\write\gtoutfile{}
\ifx\theasciiabstract\relax
\immediate\write\gtoutfile{\theabstract}\else
\immediate\write\gtoutfile{\theasciiabstract}\fi
\immediate\write\gtoutfile{}
\immediate\write\gtoutfile{\noexpand\\}
\immediate\write\gtoutfile{}
\immediate\closeout\gtoutfile}}  %%% end of definition of \makeheadfile

\def\maketitlepage{\makeagttitle\makeheadfile}
\let\makeshorttitle\maketitlepage
\let\maketitle\maketitlepage

\def\ifplaintex{\expandafter\ifx\csname documentclass\endcsname\relax}

\def\gtp{{\mathsurround=0pt\it $\cal G\mskip-2mu$eometry \&\ 
$\cal T\!\!$opology $\cal P\!$ublications}}  % GT publications

\def\recd{{\small Received:\qua\receiveddate\ifx\reviseddate\relax
\else\qquad Revised:\qua\reviseddate\fi\par}} 

\def\lognumber#1{\def\thelognumber{#1}}
\def\volumenumber#1{\def\thevolumenumber{#1}}
\def\volumeyear#1{\def\thevolumeyear{#1}}
\def\papernumber#1{\def\thepapernumber{#1}}
\def\pagenumbers#1#2{\def\startpage{#1}\def\finishpage{#2}}
\def\published#1{\def\publishdate{#1}}

\def\received#1{\def\receiveddate{#1}}
\def\revised#1{\def\reviseddate{#1}}
\def\accepted#1{\def\accepteddate{#1}}

\def\asciiaddress#1{\def\theasciiaddress{#1}}

\long\def\asciiabstract#1{\long\def\theasciiabstract{#1}}

\let\\\par\let\thelognumber\relax\let\thevolumenumber\relax
\let\thepapernumber\relax\let\thevolumeyear\relax\let\startpage\relax
\let\finishpage\relax\let\publishdate\relax\let\receiveddate\relax
\let\reviseddate\relax\let\accepteddate\relax\let\theasciititle\relax
\let\theasciiauthors\relax\let\theasciiaddress\relax
\let\theasciiabstract\relax

\let\theasciiemail\relax

\ifplaintex
\font\logobig=cmssbx10 scaled 3836
\font\logomed=cmssbx10 scaled 2557
\else
\font\logobig=cmssbx10 scaled 4200
\font\logomed=cmssbx10 scaled 2800
\fi

\long\def\makeagttitle{   %%% start of definition of \makeagttitle
\count0=\startpage
\agt\hfill      %   Journal title (top left) 
\hbox to 45truept{\vbox to 0pt{\vglue -13truept{\logomed A\kern -.37em{\logobig 
T}\kern -.38em G}\vss}\hss}
\break
{\small Volume \thevolumenumber\ (\thevolumeyear)
\startpage--\finishpage\nl
Published: \publishdate}

\vglue .25truein

{\parskip=0pt\leftskip 0pt plus
1fil\def\\{\par\smallskip}{\Large\bf\thetitle}\par\medskip} \vglue
0.05truein

{\parskip=0pt\leftskip 0pt plus 1fil\def\\{\par}{\sc\theauthors}
\par\medskip}%
 
\vglue 0.03truein 

{\small\leftskip 25truept\rightskip 25truept{\bf Abstract}\stdspace\theabstract

{\bf AMS Classification}\stdspace\theprimaryclass
\ifx\thesecondaryclass\relax\else; \thesecondaryclass\fi\par
{\bf Keywords}\stdspace \thekeywords\par}\vglue 7truept

}   %%%% end of definition of \makeagttitle

\ifplaintex
\font\phead=cmsl9 scaled 950
\font\pnum=cmbx10 scaled 913
\font\pfoot=cmsl9 scaled 950
\headline{\vbox to 0pt{\vskip -4.5mm\line{\small\phead\ifnum
\count0=\startpage ISSN 1472-2739 (on-line) 1472-2747 (printed)
\hfill {\pnum\folio}\else\ifodd\count0\def\\{ }% 
\ifx\theshorttitle\relax\thetitle\else\theshorttitle\fi\hfill{\pnum\folio}
\else\def\\{ and }{\pnum\folio}\hfill\ifx\theshortauthors\relax\theauthors
\else\theshortauthors\fi\fi\fi}\vss}}
\footline{\vbox to 0pt{\vglue 0mm\line{\small\pfoot\ifnum\count0=\startpage
\copyright\ \gtp\hfill\else
\agt, Volume \thevolumenumber\ (\thevolumeyear)\hfill\fi}\vss}}
\else
\font=cmsl9 scaled 1050
\font\lnum=cmbx10 
\font=cmsl9 scaled 1050
\makeatletter
\def\@oddhead{{\small\lhead\ifnum\count0=\startpage ISSN 1472-2739 
(on-line) 1472-2747 (printed)\hfill {\lnum\number\count0}\else\ifodd\count0
\def\\{ }\ifx\theshorttitle\relax \thetitle \else\theshorttitle\fi\hfill
{\lnum\number\count0}\else\def\\{ and }{\lnum\number\count0}
\hfill\ifx\theshortauthors\relax 
\theauthors\else\theshortauthors\fi\fi\fi}}\def\@evenhead{\@oddhead}
\def\@oddfoot{\small\lfoot\ifnum\count0=\startpage\copyright\ \gtp\hfill\else
\agt, Volume \thevolumenumber\ (\thevolumeyear)\hfill\fi}
\def\@evenfoot{\@oddfoot}
\makeatother
\fi
\let\maketitlepage\makeagttitle
\let\makeshorttitle\maketitlepage
\let\maketitle\maketitlepage

\newwrite\gtoutfile
\long\gdef\makeheadfile{  %%% start of definition of \makeheadfile
{\def\\{, }\def\s{ }
\immediate\openout\gtoutfile head.xxx
\immediate\write\gtoutfile{To: math@arxiv.org}
\immediate\write\gtoutfile{Subject: put OR rep NNNNN:ppppp}
\immediate\write\gtoutfile{--text follows this line--}
\immediate\write\gtoutfile{Proxy-for: \ifx\theasciiauthors\relax
\theauthors\else\theasciiauthors\fi\s<\ifx\theasciiemail\relax\theemail\else\theasciiemail\fi>}
\immediate\write\gtoutfile{\noexpand\\}
\immediate\write\gtoutfile{Authors: \ifx\theasciiauthors\relax
\theauthors\else\theasciiauthors\fi}
{\def\\{ }\immediate\write\gtoutfile{Title: \ifx\theasciititle\relax
\thetitle\else\theasciititle\fi}}
\immediate\write\gtoutfile{Subj-class: GT or SG, GR etc}
\immediate\write\gtoutfile{MSC-class: \theprimaryclass\ifx\thesecondaryclass\relax\else, \thesecondaryclass\fi}
\immediate\write\gtoutfile{Journal-ref: Algebr. Geom. Topol. \thevolumenumber\s
(\thevolumeyear) \startpage-\finishpage}
\immediate\write\gtoutfile{Comments: Published by Algebraic and
Geometric Topology at}
\immediate\write\gtoutfile{\s\s\s  http://www.maths.warwick.ac.uk/agt/AGTVol\thevolumenumber/agt-\thevolumenumber-\thepapernumber.abs.html}
\immediate\write\gtoutfile{\noexpand\\}
\immediate\write\gtoutfile{}
\ifx\theasciiabstract\relax
\immediate\write\gtoutfile{\theabstract}\else
\immediate\write\gtoutfile{\theasciiabstract}\fi
\immediate\write\gtoutfile{}
\immediate\write\gtoutfile{\noexpand\\}
\immediate\write\gtoutfile{}
\immediate\closeout\gtoutfile}}  %%% end of definition of \makeheadfile

\def\maketitlepage{\makeagttitle\makeheadfile}
\let\makeshorttitle\maketitlepage
\let\maketitle\maketitlepage

\lognumber{37}
\volumenumber{2}
\volumeyear{2002}
\papernumber{37}
\pagenumbers{921}{936}
\received{6 May 2002}
\revised{16 September 2002}
\accepted{12 October 2002}
\published{21 October 2002}

\usepackage{amssymb, graphicx, amsmath}    %% newams

\newtheorem{thm}{Theorem}[section]  
\newtheorem{lem}[thm]{Lemma}        
\newtheorem{prop}[thm]{Proposition} 
\newtheorem{corollary}[thm]{Corollary}  
\theoremstyle{definition}

\def\Cal#1{{\cal#1}}
\def\<{\langle}\def\>{\rangle}\def\what{\widehat}\def\wtil{\widetilde}
\def\Z{{\mathbb Z}}\def\N{{\mathbb N}} 
\def\R{{\mathbb R}} \def\E{{\mathbb E}}

\let\Proof\proof
\def\Remark{\par\medskip{\bf Remark}\qua}

\def\vert{{\cal V}}\def\edge{{\cal E}}\def\gen{{\cal G}}
\def\dss{\text{dim}_{ss}}\def\cd{\text{cd}}\def\gd{\text{gd}}
\def\Min{\text{\sl Min}}

\def\al{\alpha}                 \def\be{\beta}
			\def\Ga{\Gamma}
\def\sig{\sigma}

\begin{document}

\title{On the CAT(0) dimension of 2-dimensional\\Bestvina-Brady groups}                    
\authors{John Crisp}

\address{Laboratoire de Topologie, Universit\'e de Bourgogne\\
UMR 5584 du CNRS -- BP 47 870, 21078 Dijon, France}                  
\asciiaddress{Laboratoire de Topologie, Universite de Bourgogne\\UMR 
5584 du CNRS - BP 47 870, 21078 Dijon, France}                  

\email{jcrisp@u-bourgogne.fr}

\begin{abstract}  
Let $K$ be a $2$-dimensional finite flag complex. We study the CAT(0) dimension
of the `Bestvina-Brady group', or `Artin kernel', $\Ga_K$. We show that $\Ga_K$ has 
CAT(0) dimension $3$ unless $K$ admits a piecewise Euclidean metric of 
non-positive curvature. We give an example to show that this implication cannot
be reversed. Different choices of $K$ lead to examples where the CAT(0) dimension 
is $3$, and either (i) the geometric dimension is $2$, or (ii) the cohomological
dimension is $2$ and the geometric dimension is not known.
\end{abstract}

\asciiabstract{Let K be a 2-dimensional finite flag complex. We study
the CAT(0) dimension of the `Bestvina-Brady group', or `Artin kernel',
Gamma_K.  We show that Gamma_K has CAT(0) dimension 3 unless K admits a
piecewise Euclidean metric of non-positive curvature. We give an
example to show that this implication cannot be reversed. Different
choices of K lead to examples where the CAT(0) dimension is 3, and
either (i) the geometric dimension is 2, or (ii) the cohomological
dimension is 2 and the geometric dimension is not known.}


\maketitle

Let $\Gamma$ be a countable group. We denote the cohomological dimension of $\Gamma$ by
$\cd(\Gamma)$ and the geometric dimension by $\gd(\Gamma)$. See \cite{Bro}, pp 184--5,
for definitions of these notions. It was shown by Eilenberg and Ganea \cite{EG} that if
$\cd(\Gamma)\geq 3$ then $\cd(\Gamma)=\gd(\Gamma)$. The same is true if $\cd(\Gamma)=0$
($\Gamma$ trivial), or $\cd(\Gamma)=1$ (by the Stallings--Swan Theorem \cite{St,Sw}).
There remains however the possibility that there exists a group $\Gamma$ with
$\cd(\Gamma)=2$ but $\gd(\Gamma)=3$, that is, a counter-example to the so-called
Eilenberg-Ganea Conjecture. In the search for such an example, two promising
families of groups have been proposed. Both constructions begin by choosing a finite
flag $2$-complex $K$ which is acyclic but not contractible, such as, for example,
any finite flag subdivision of the 2-spine of the Poincar\'e homology sphere.
For the first family, due to Bestvina and Davis, one considers a torsion free finite
index subgroup $G_K$ of the right-angled  Coxeter group associated to the $1$-skeleton
$K^1$. In the second case (Bestvina-Brady \cite{BB}) one considers the kernel
$\Gamma_K$ of the length homomorphism on the right-angled Artin group associated
to $K^1$. These groups may be defined for an arbitrary finite flag complex $K$.
In the case of interest, where $K$ is a non-contractible, acyclic flag $2$-complex,
the groups $G_K$ and $\Gamma_K$ are known to  have cohomological dimension $2$ 
but their geometric dimensions are not known. 

Metric spaces of global nonpositive curvature, or CAT(0) spaces, provide a  natural
supply of contractible spaces. Moreover, in many cases the CAT(0) condition, and hence
contractibility of a given space, may be easily verified in terms of local geometric information.
This motivates our interest in CAT(0) spaces and group actions on these spaces, and leads
us to the notion of the \emph{CAT(0) dimension} of a group
(we adopt the definition used by Bridson in \cite{Bri2}):
$$\phantom{\rm where}\qquad\dss(\Gamma)=min(\Cal D\cup\{\infty\})\,,\qquad {\rm where}$$
\begin{gather*}
\Cal D= \{ \dim(X) : \text{$X$ is a } \text{complete CAT(0) space on which}\\
\phantom{\rm where}\hspace{1in}
\Gamma\text{ acts properly by semi-simple isometries}\}\,.
\end{gather*}
Here $\dim(X)$ denotes the covering dimension of a metric space $X$ (see \cite{HW}).
Note that $\Cal D$ may well be an empty set (i.e.\ no 
such actions exist) in which case $\dss(\Gamma)=\infty$.

Let $K$ be an arbitrary finite $2$-dimensional flag complex. The main result
of this paper (Theorem \ref{main}) states that, unless the $2$-complex $K$ is
aspherical (in fact, unless it admits a piecewise Euclidean metric of nonpositive
curvature), the group $\Gamma_K$ does not act properly semi-simply on any
$2$-dimensional CAT(0) space, and hence has CAT(0) dimension $3$.
In particular, $\dss(\Gamma_K)=3$ in the case where $K$ is
acyclic but not aspherical, as for example when $K$ is a flag decomposition of 
the 2-spine of the Poincar\'e homology sphere (Corollary \ref{cor5}). However, we
do not rule out the possibility that $\dss(\Gamma_K)=2$ for some aspherical but
non-contractible flag $2$-complex $K$.
We note also that this result is still a very long way from showing that some $\Gamma_K$
with cohomological dimension $2$ has geometric dimension $3$. It merely indicates,
as one might already expect, that for these examples a $2$-dimensional Eilenberg-MacLane
complex will not be so easy to find, if in fact such a complex exists.

Theorem \ref{main} is a result of the fact that the group $\Gamma_K$ contains many
`overlapping' $\Z\times\Z$ subgroups, which imply the existence of flat planes
embedded in any CAT(0) space on which the group acts. Our argument proceeds by
studying how these flat planes interact.  
The Bestvina-Davis examples, $G_K$, also contain $\Z\times\Z$ subgroups, however
they arise in a less regular manner than in the Artin kernels, and the present techniques 
seem less applicable. It would be interesting to know whether similar results
hold for the groups $G_K$.

It is already known that requiring that a group act properly semi-simply on a CAT(0)
space imposes constraints both on the group and on the dimension of the space.
For example, the Baumslag-Solitar groups
$BS(n,m)=\langle x,t | t^{-1}x^nt=x^m\rangle$, for $1\leq n<m$, have geometric
dimension $2$ but do not admit any proper semi-simple actions on CAT(0) spaces
(which may be seen by considering the translation length of $x$).
There are also recent examples given of finitely presented groups $\Gamma$ having
$\gd(\Gamma)=2$ but $\dss(\Gamma)=3$ (see \cite{BC,Bri2}).
We note that Bridson's example
\cite{Bri2} is actually an index two subgroup of a group with CAT(0) dimension $2$.
The examples of both \cite{BC} and \cite{Bri2} are all \emph{CAT(0) groups}, that is
they act properly and cocompactly on CAT(0) spaces.

In Corollary \ref{cor2} we give a simple method for constructing infinitely many further
examples of finitely presented groups with geometric dimension $2$ but CAT(0) dimension
$3$. These are the groups $\Gamma_K$ where $K$ is any flag triangulation of a
contractible $2$-complex which does not admit a CAT(0) metric. 
Zeeman's dunce hat is such a complex.
We note that the 3-dimensional CAT(0) actions which are known to exist for these
groups are semi-simple, but not cocompact. It is not known whether any of the
examples given by Corollary \ref{cor2} are CAT(0) groups. 

In the last section of the paper we give a refinement (Theorem \ref{refined}) of the
main result which, in the case that $K$ is a simply-connected finite flag $2$-complex,
leads to an ``if and only if'' statement. This allows us to give an example of a
group $\Gamma_K$  where $K$  admits a CAT(0) piecewise Euclidean metric,
but where we still have $\gd(\Gamma_K)=2$ and $\dss(\Gamma_K)=3$.
Since submission of this paper, Noel Brady has informed
me that this particular example has a cubic Dehn function \cite{Bra},
and so cannot be a CAT(0) group.

We close with the following remark: all the groups $\Gamma$ considered in 
this paper have the property that every finite index subgroup of $\Gamma$ contains a subgroup 
isomorphic to $\Gamma$. Thus all the examples given in this paper 
of 2-dimensional groups with CAT(0) dimension 3 have the further property
that every finite index subgroup has CAT(0) dimension 3 as well. 
This distinguishes our examples from the example given by Bridson in [4] 
which has an index 2 subgroup with CAT(0) dimension 2.
On the other hand, it is not known whether the CAT(0) dimension can drop
when passing to finite index subgroups 
in any of the Artin group examples given in [2].

\medskip
{\bf Acknowledgement}\qua The question addressed in this paper arose
out of a talk given by Noel Brady in Dijon, June 2000. I would 
like to thank Noel Brady, Jens Harlander, Ian Leary and the referee for many helpful
comments and contributions.
This paper was largely written in 2001 during a 6 month stay at the UMPA ENS-Lyon.
I acknowledge the support of the CNRS during that period, and also wish to
thank the members of the Unit\'e de Math\'ematiques Pures et Appliqu\'ees, ENS-Lyon,
for their hospitality.

\section{Some CAT(0) geometry}\label{Sect1}

We introduce just the basic concepts of CAT(0) geometry which we will need here.
We refer to \cite{BH} for details.
 
Let $(X,d)$ be a metric space. By a \emph{geodesic segment} with endpoints $x$ and $y$,
usually denoted $[x,y]$, we mean the image of an isometric embedding of the closed
interval $[0,d(x,y)]$ into $X$ such that the endpoints are mapped to $x$ and $y$. 
A \emph{geodesic triangle} $\Delta(x,y,z)$ in  $X$ 
is simply a union of three geodesic segments $[x,y] \cup [y,z] \cup [z,x]$. 
Note that $\Delta$ need not be uniquely determined by the three points $x,y$ and $z$.
However, for every triple $x,y,z\in X$ there is (up to isometry) a unique
\emph{comparison triangle} $\Delta'(x,y,z)$ in the Euclidean plane $\E^2$
which has vertices $x',y',z'$ such that $d_{\E^2}(x',y')=d(x,y)$, $d_{\E^2}(y',z')=d(y,z)$,
and $d_{\E^2}(x',z')=d(x,z)$. 
In the case that $x$ is distinct from both $y$ and $z$, we write $\angle'_x(y,z)$ for the angle
at $x'$ between the sides $[x',y']$ and $[x',z']$ in the comparison triangle $\Delta'$.

There are various equivalent formulations of global nonpositive curvature (for geodesic spaces)
in terms of comparison triangles (see \cite{BH}, Chapter II.1). We will make use of the following one:

Given distinct points $x,y,z\in X$ and geodesic segments $[x,y],[x,z]$ one defines the 
\emph{(Alexandrov) angle} $\angle([x,y],[x,z])$ between the two geodesic segments to be
\[
\lim_{\epsilon\to 0}\ \sup \{ \angle'_x(p,q) : p\in[x,y],q\in[x,z]
\text{ and } 0<d(p,x),d(q,x)<\epsilon \}\ .
\]
A geodesic metric space $X$ (one in which there is a geodesic segment joining any pair
of points) is said to be a \emph{CAT(0) space} if, for every geodesic triangle
$\Delta(x,y,z)$ with distinct vertices, the Alexandrov angle between any two sides of
$\Delta$ is no greater than the corresponding angle in the comparison triangle
$\Delta'(x,y,z)$ in $\E^2$.

This definition leads quickly to the following:

\begin{lem}\label{polygon}
Let $X$ be a CAT(0) space, and $P=P(v_1,v_2,\ldots ,v_m)$ a geodesic polygon in $X$
with distinct vertices $v_1,v_2,\ldots ,v_m$ ($m\geq 3$) and sides $[v_i,v_{i+1}]$
for $1\leq i\leq m$ (indices mod $m$). Write $\phi_i=\angle([v_i,v_{i-1}],[v_i,v_{i+1}])$ 
for the angle of the polygon at $v_i$, for $i=1,\ldots ,m$. Then 
\[
\text{angle sum of } P\  =\  \sum\limits_{i=1}^m
\phi_i\ \  \leq \ (m-2)\pi \,.
\]
\end{lem}

\Proof
 If $m=3$ then the result follows directly from the CAT(0) inequality. Otherwise,
take the (unique) geodesic segment $[v_1,v_3]$ to form a triangle $\Delta(v_1,v_2,v_3)$
and an $(m-1)$-gon $P'(v_1,v_3,\ldots ,v_m)$. Let $\psi_1$, $\psi_2(=\phi_2)$ and $\psi_3$
be the angles of $\Delta$, and $\phi'_1,\phi'_3,\phi_4,\ldots ,\phi_m$ the angles of $P'$.
Then, by induction on the number of sides of a polygon, we have 
\begin{gather}\label{eq1}
\phantom{\rm and}\qquad \psi_1+\phi_2+\psi_3\leq\pi\,,\qquad {\rm and}\\
\label{eq2}\phi'_1+\phi'_3+\phi_4+ \cdots +\phi_m\leq (m-3)\pi
\end{gather}
But $\phi_1\leq\phi'_1+\psi_1$ and $\phi_3\leq\phi'_3+\psi_3$, so that the result follows
immediately by adding the two inequalities (\ref{eq1}) and (\ref{eq2}).\endproof

 Let $g$ be an isometry of a 
metric space $(X,d)$. The {\it translation length of $g$}, 
denoted by $\ell(g)$, is defined to be 
$$
\ell(g) \; = \; \inf\{d(x,gx) \, | \, x \in X\} \, ,
$$
and the {\it minset of $g$}, denoted by $\Min(g)$, is defined to be the 
possibly empty set 
$$
\Min(g) \; = \; \{x \in X \, | \, d(x,gx) = \ell(g)\}\,. 
$$
An isometry of a CAT(0) space is called {\it semi-simple} 
if it has a nonempty minset. If $g$ is semi-simple with nonzero translation length
then we say that $g$ is \emph{hyperbolic}. In this case any $g$-invariant geodesic
line in $X$ shall be called an \emph{axis of $g$} or \emph{$g$-axis}.
By \cite{BH} Proposition II.6.2, the $g$-axes all lie in $\Min(g)$.

The following statements may be found in \cite{BH}, Chapter II.6. Let $g$ be a
semi-simple isometry of a CAT(0) space $X$. Then $\Min(g)$ is always a closed, convex
subspace of $X$. If $\ell(g)\neq 0$ then $\Min(g)$ is isometric to the metric product of
$\R$ with a CAT(0) space $Y$, where each $\R$-fibre is an axis of $g$. The element $g$
acts on $\Min(g)$ by translating along the ${\R}$ factor and fixing the $Y$ factor
pointwise, and each isometry of $X$ which commutes with $g$ leaves $\Min(g)$ invariant
while respecting the product structure.
Moreover, if $g$ belongs to a group $\Gamma$ of semi-simple isometries acting properly
on $X$ then $C_\Gamma(g)/\<g\>$ acts properly by semi-simple isometries on $Y$.
(c.f: Propositions II.6.9-10 of \cite{BH}). Here $C_\Gamma(g)$ denotes the
centralizer of $g$ in $\Gamma$.

We wish to restrict our attention now to group actions on $2$-dimensional spaces.
The following is a special case of Bridson's Flat Torus Theorem \cite{Bri1,BH}.
We refer to \cite{HW} for the theory of `covering dimension'. We note that the
covering dimension of a simplicial complex agrees with the (simplicial) dimension,
and that the covering dimension of a metric space is bounded below by the dimension of any subspace.

\begin{prop}\label{flat}
Let $A$ be a free abelian group of rank $2$ acting properly by semi-simple isometries
on a CAT(0) space $X$ of covering dimension $2$.  Then $\Min(A) = \cap_{a \in A}\Min(a)$
is an $A$-invariant isometrically embedded flat plane ($\cong\E^2$) and the group $A$
acts by translations on $\Min(A)$ with quotient a $2$-torus.
\end{prop}

\Proof 
By the Flat Torus Theorem, $\Min(A)$ is non-empty and splits as a product $Y\times\E^2$,
where each fibre $\{y\}\times\E^2$ has the desired properties. The dimension constraint
ensures that $Y$ consists of precisely one point. For if $p,q\in Y$ are distinct points
then, by convexity, $Y$ contains the geodesic segment $[p,q]$, and hence $X$ contains a
subspace $[p,q]\times\E^2$ of dimension $3$, contradicting $\text{dim}(X)=2$.
\endproof

\begin{notation}
Suppose that $A$ is a free abelian group of rank $2$ acting properly 
by semi-simple isometries on a CAT(0) space $X$ of covering dimension $2$.
For any pair of nontrivial elements $x,y\in A$, we write $\theta(x,y)$
for the angle between any positively oriented $x$-axis and any positively
oriented $y$-axis in $\Min(A)$. More precisely
\[
\theta(x,y)=\angle([p,x(p)],[p,y(p)])\,, \text{ for any } p\in \Min(A)\,. 
\]
Note that $0<\theta(x,y)<\pi$ in the case that $x$ and $y$
generate $A$.
\end{notation}

We now prove the following Lemma.

\begin{lem}\label{thetas}
Suppose that $\Gamma$  is a group acting properly by semi-simple
isometries on a CAT(0) space $X$ of covering dimension $2$. 
Suppose that we have a cyclic sequence of $m \geq 4$ group elements
$a_1,a_2,\ldots ,a_m=a_0$ of $\Gamma$ (indices taken mod $m$) such that,
for each $i=1,\ldots ,m$, we have

\begin{itemize}
\item[\rm(i)] $\<a_{i-1},a_{i+1}\>= F_2$ (free group of rank 2), and
\item[\rm(ii)] $\< a_i, a_{i+1}\>\cong \Z\times\Z$.
\end{itemize}

Then\qquad\qquad\qquad\qquad\qquad\qquad
$\displaystyle
\sum\limits_{i=1}^m \theta(a_i,a_{i+1}) \geq 2\pi\,.
$
\end{lem}

\Proof
For $i=1,2,\ldots ,m$, write $\Pi(i,i+1)$ for the flat plane 
$\Min(\< a_i,a_{i+1}\>)$ in $X$ (see Proposition \ref{flat}). 
Let $C(i)$ denote the convex closure of $\Pi(i-1,i)\cup\Pi(i,i+1)$ in $X$.

Fix $i\in\{ 1,\ldots ,m\}$. 
Note that $C(i)$ is contained in $\Min(a_i)$ (by convexity of the minset), and that
$\Min(a_i)$ is isometric to $T\times\R$ where, since $X$ is $2$-dimensional, $T$ must
be an $\R$-tree (see Lemma 3.2 of \cite{Bri2} for more detail). With respect to this
decomposition, $\Pi(i-1,i)\cong\tau\times\R$ and $\Pi(i,i+1)\cong\sig\times\R$, where
$\tau$ and $\sig$ are the axes in $T$ for the elements $a_{i-1}$ and $a_{i+1}$
respectively, and  $C(i)\cong H\times\R$ where $H$ is the convex closure in $T$
of $\tau\cup\sig$.

By hypotheses (i) and (ii) the elements $a_{i-1}$ and $a_{i+1}$ generate a rank $2$
free group which commutes with $a_i$ and therefore acts properly semi-simply on $T$.
It follows that their axes $\tau$ and $\sig$ in $T$ are either disjoint, or intersect
in a closed interval of finite length. We define the finite (possibly zero) length
closed interval $I=I(i)$ in $T$ as follows. If $\tau$ and $\sig$ are disjoint,
then $I$ is the unique shortest geodesic segment joining them. In this case
$H=\tau\cup I\cup\sig$. Otherwise, we set $I=\tau\cap\sig$.

Now define the subspace $A(i)=I\times\R$ of $C(i)$, noting that $A(i)$ separates each
of $\Pi(i-1,i)$ and $\Pi(i,i+1)$ into two  $\< a_i\>$-invariant components. 

For each $i=1,\ldots ,m$, define $Q(i,i+1)$ to be the 
component of $\Pi(i,i+1)\setminus(A(i)\cup A(i+1))$ such that both 
$a_i(Q(i,i+1))$ and $a_{i+1}(Q(i,i+1))$ lie again in  $Q(i,i+1)$.
This is a sector of the plane with angle $\theta(a_i,a_{i+1})$ and bounded 
by positive semi-axes for $a_i$ and $a_{i+1}$ respectively.

Again fix $i\in\{ 1,\ldots ,m\}$. Now, for some $p,q\in I(i)$ and $r,s\in\R$, we may write
\[
\begin{aligned}
Q(i-1,i)\cap A(i) &= \{p\}\times\{ t\in\R : t\geq r\}\,,\text{ and}\\
Q(i,i+1)\cap A(i) &= \{q\}\times\{ t\in\R : t\geq s\}\,.
\end{aligned}
\]  
Let $t_0=\text{max}(r,s)$ and define 
\[
\begin{aligned}
B&(i)=[p,q]\times \{ t\in\R : t\geq t_0\}\,,\\
\mu&(i)=\{p\}\times \{ t\in\R : t\geq t_0\}\,,\ b_i=(p,t_0)\in\mu(i)\,, \text{ and}\\
\what Q&(i,i+1)=B(i)\cup Q(i,i+1)\,.
\end{aligned}
\] 
These definitions are illustrated in Figure \ref{fig0} below.
Note that the ray $\mu(i)$ is common to both $\what Q(i-1,i)$ and $\what Q(i,i+1)$.
Also, $\what Q(i,i+1)$ is contained in the convex set $C(i)$.We now choose, for
each $i=1,\ldots ,m$, a point $v_i$ in $\mu(i)$, different from the basepoint $b_i$
of the ray $\mu(i)$. These points may be chosen sufficiently far along their
corresponding rays that the unique geodesic segment in $X$ from $v_i$ to $v_{i+1}$ 
lies wholly in $\what Q(i,i+1)$. It is also easy to arrange that the $v_i$ are
mutually distinct.
We now apply Lemma \ref{polygon} to the polygon
$P=P(v_1,v_2,\ldots ,v_m)$. See Figure \ref{fig0}.

\begin{figure}[ht]
\begin{center}
\includegraphics[width=13cm]{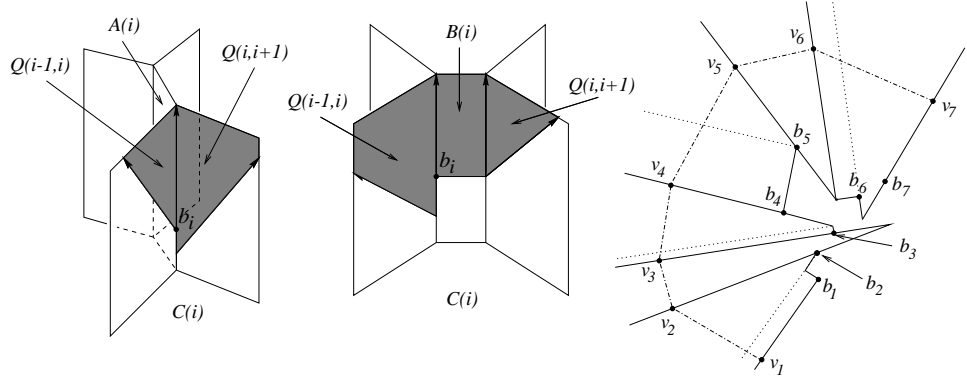}
\end{center}
\vglue -10pt
\caption{Examples of $C(i)$ and the construction of the polygon $P(v_1,..,v_m)$}\label{fig0}
\end{figure}

For each $i=1,\ldots ,m$, write $$\al_i=\angle([v_i,v_{i-1}],[v_i,b_i])
\text{\ \ \ and \ \ }\be_i=\angle([v_i,v_{i+1}],[v_i,b_i]).$$ From the Euclidean geometry of $\Pi(i,i+1)$ we
have, for each $i$, 
\begin{equation}\label{one}
\be_i+\al_{i+1}+\theta(i,i+1)=\pi\,.
\end{equation}
Note that by choosing the points $v_i$ equidistant from their
corresponding basepoints $b_i$ and sufficiently far away, we may suppose that the angles $\al_i$ and
$\be_i$ are all acute (in fact $\alpha_{i+1}$ and $\be_i$ approach one another as
$d(v_i,b_i)=C\to\infty$, since for any fixed point $b$ in the vicinity of either $b_i$ or $b_{i+1}$,
the triangle $\Delta(v_i,v_{i+1},b)$ will tend towards an isosceles triangle). 

As in Lemma  \ref{polygon}, write $\phi_i=\angle([v_i,v_{i-1}],[v_i,v_{i+1}])$.
By measuring in the convex subspace $C(i)$ of $X$, and supposing that $\al_i$ and $\be_i$ are acute,
it is easily seen that 
\begin{equation}\label{two}
\phi_i = \al_i +\be_i\,.
\end{equation}
Combining (\ref{one}) and (\ref{two}), and applying Lemma \ref{polygon}, we now have
\[
\sum\limits_{i=1}^m \theta(i,i+1)\  =\  m\pi - \sum\limits_{i=1}^m\phi_i\ 
\geq\  2\pi\,,
\]
completing the proof.
\endproof

\Remark We have tried to make the proof of the above Lemma as elementary as possible. 
The argument may be simplified considerably by using the fact that, if $X$ is
a complete CAT(0) space, then its boundary $\partial X$ (with the Tits, or angular,
metric) is a complete CAT(1) space (see Theorem II.9.13 of \cite{BH}).
Each flat sector $Q(i,i+1)$ determines a geodesic
segment of length $\theta(i,i+1)$ between the points $a_i^{+\infty}$ and
$a_{i+1}^{+\infty}$ in $\partial X$. The essence of the argument is to see that the
union of these segments forms a locally geodesic circle in $\partial X$ and hence
must have total length at least $2\pi$. (Local geodicity at $a_i^{+\infty}$ comes from
the fact that, by the dimension constraint, the minset $\Min(a_i)$ is isometric to
$T\times \R$, where $T$ is an $\R$-tree, and is convex in $X$).
Note however that the more elementary proof
makes no use at all of the hypothesis that $X$ be complete.

\section{Application to Artin kernels}

A simplicial complex is called a \emph{flag complex} if whenever $n+1$ vertices span a
complete graph in the 1-skeleton they span an $n$-simplex in $K$.
Let $K$ be a finite connected $2$-dimensional flag complex. Let $\vert(K)$ denote the
set of vertices of $K$ and $\edge(K)$ the set of ordered pairs
$(u,v)\in\vert(K)\times\vert(K)$ such that $u$ and $v$ are distinct vertices joined by
an edge in $K$. The elements $(u,v)\in\edge(K)$ are thought of as directed edges of $K$.

The \emph{right-angled Artin group} associated to the $1$-skeleton of $K$ is the group
\[
A_K = \< \vert(K)\mid uv=vu \text{ if } (u,v)\in\edge(K)\>\,.
\]
The fact that $K$ is a $2$-dimensional flag complex ensures that mutually commuting
subsets of the generating set $\vert(K)$ contain at most three elements, and that whenever
one sees three distinct mutually commuting generators they span a $2$-simplex in $K$.
In fact $A_K$ acts freely on a $3$-dimensional CAT(0) cube complex \cite{BB}, so that 
\[
\cd(A_K)=\gd(A_K)=\dss(A_K)=3\, .
\]
Let $l:A_K\to\Z$ denote the `length' homomorphism which takes every generator
in $\vert(K)$ to $1$. We define the \emph{Artin kernel}, or \emph{Bestvina-Brady group}, to be the group
\[
\Gamma_K=\ker(l:A_K\to\Z)\,.
\]
Since $K$ is a connected complex, it follows (see \cite{DL}) that $\Gamma_K$
is generated by the set
\[
\gen(K)=\{ x_{(u,v)}=u^{-1}v : (u,v)\in\edge(K) \}\,.
\]
In fact, in \cite{DL}, Dicks and Leary show that a complete presentation of $\Gamma_K$
may be obtained by taking as relators all words $x^n_{(u_1,u_2)}x^n_{(u_2,u_3)}\cdots x^n_{(u_k,u_1)}$
where the sequence of edges $(u_1,u_2),(u_2,u_3),\ldots ,(u_k,u_1)\in \edge(K)$ forms a
directed cycle in $K$, and where $n\in\Z$. In particular, we have the following relations in $\Gamma_K$:
\[
\begin{aligned}
x_{(v,u)}=x_{(u,v)}^{-1}\hskip5mm &\text{for all } (u,v)\in\edge(K)\,,\hfil\\
x_{(u,v)}x_{(v,w)}=x_{(u,w)}=x_{(v,w)}x_{(u,v)}
\hskip5mm &\text{whenever $u,v,w$ span a $2$-simplex in $K$}\,.
\end{aligned}
\]
It is known, from \cite{vdL}, that if $L$ is a full subcomplex of $K$, then the natural map $A_L\to A_K$
is injective. For $L<K$, we identify $\Gamma_L$ and $A_L$  with their images in $\Gamma_K$ and $A_K$
respectively, via this map. In particular, if $\sigma$ is a $2$-simplex in $K$, then $\Gamma_\sigma$ is
a free abelian group of rank $2$ sitting inside $A_\sigma\cong\Z\times\Z\times\Z$. 

Before stating the following Theorem we recall that a geodesic metric space is said to have
\emph{nonpositive curvature} if it is locally CAT(0), i.e: if every point has a convex open
neighbourhood which is CAT(0) with the induced metric. If $X$ is a complete geodesic metric
space of nonpositive curvature then its universal cover $\wtil X$ is a CAT(0) space (with the
induced length metric). This is a generalisation of the Cartan-Hadamard Theorem (see \cite{BH}
Chapter II.4 for a discussion). Any CAT(0) space is contractible (\cite{BH} II.1.5), so that in
the above situation $\wtil X$ is contractible and $X$ is aspherical.  We now have the following:

\begin{thm}\label{main}
Let $K$ be a finite connected $2$-dimensional flag complex.
If $\Gamma_K$ acts properly by semi-simple isometries on a $2$-dimensional 
CAT(0) space $X$, then $K$ admits a complete piecewise Euclidean metric of 
nonpositive curvature, and in particular, $K$ is aspherical. 
\end{thm}

\Proof
We define a piecewise Euclidean metric on $K$ as follows. Firstly define
a Euclidean metric on each $2$-simplex $\sigma$ of $K$ in such a way that the length
of each edge $(u,v)$ of $\sigma$ is equal to the translation length of $x_{(u,v)}$ 
on $X$. (That this is always possible will become clear in the following paragraph).
Now, as in \cite{BH} Chapter I.7, there is a complete geodesic metric on $K$ defined
by setting $d(x,y)$ to be the infimum of the lengths of all piecewise linear paths
from $x$ to $y$ in $K$ (where the length of such a path is just the sum of the
lengths of its linear segments as measured inside the individual simplexes). 

Suppose that $u,v,w$ span a $2$-simplex $\sig$ in $K$. The elements 
$x_{(u,v)}$, $x_{(u,w)}$ and $x_{(v,w)}$ generate a rank two abelian subgroup
$\Gamma_\sigma$ of $\Gamma_K$ so that, by Proposition \ref{flat}, their minsets 
intersect in a $\Gamma_\sigma$-invariant flat plane $\Pi$.
Take any point $W\in\Pi$ and let $\Delta$ denote the triangle in $\Pi$ with
vertices $W$, $V=x_{(v,w)}(W)$ and $U=x_{(u,w)}(W)$.
Then by the `triangle relation' in $\Gamma(K)$ we have $x_{(u,v)}(V)=U$,
and clearly the triangle $\Delta(U,V,W)$ is isometric to $\sig(u,v,w)$.
In particular, $\sig$ will have an angle at $w$ precisely equal to
$\theta(x_{(v,w)},x_{(u,w)})$, and similarly for the other angles. 

We recall (\cite{BH} Chapter II.5) that a 2-dimensional piecewise Euclidean
metric complex $K$ is nonpositively curved if and only if it satisfies the
\emph{link condition}, that every simple loop in ${\rm Lk}(v,K)$ has length at
least $2\pi$, for every vertex $v$ in $K$. 
With the information given in the previous paragraph, it now follows immediately from
Lemma \ref{thetas} that the given piecewise Euclidean structure on $K$ satisfies the
link condition at every vertex, and hence is nonpositively curved. (Since $K$ is a
flag $2$-complex, the link of any of its vertices is a graph in which all circuits
have edge length at least $4$. With regard to
hypotheses (i) and (ii) of Lemma \ref{thetas}, it is not too hard to see that given
distinct edges $(u,v)$ and $(u',v')$ the elements $x_{(u,v)}$ and $x_{(u',v')}$ always
generate either $\Z\times\Z$ or a free subgroup $\Z\star\Z$ according
to whether or not the two edges lie in a common simplex).
\endproof

Let $K$ be a finite acyclic flag $2$-complex with nontrivial fundamental group.
It is known that $\cd(\Gamma_K)=2$ (Bestvina and Brady \cite{BB}), and that
$\gd(\Gamma_K)\leq\gd(A_K)=3$, but it is unknown whether or not
the geometric dimension of $\Gamma_K$ agrees with the cohomological dimension.
(Note, however, that if $K$ is acyclic and has trivial fundamental group then it is
contractible and the Bestvina-Brady argument shows that $\cd(\Gamma_K)=\gd(\Gamma_K)=2$.)
There are examples of acyclic $2$-complexes which are not aspherical. This happens, for instance,
if the acyclic complex has nontrivial torsion in its fundamental group.
One such example is the spine of the Poincar\'e homology
sphere, namely the quotient of the dodecahedral tiling of $S^2$ in which opposite faces are
identified with a $\frac{\pi}{5}$ twist. The fundamental group of this complex is finite.
Theorem \ref{main} now gives:

\begin{corollary}\label{cor5}
Let $K$ be a finite acyclic flag $2$-complex whose fundamental group has nontrivial torsion.
Then $\dss(\Gamma_K)=3$ while $\cd(\Gamma_K)=2$. In these cases $\gd(\Gamma_K)$ is not known.
\end{corollary}

If $K$ is a contractible complex then the Morse theory argument used
in \cite{BB} shows that $\gd(\Gamma_K)=2$.
In this case we have the following:

\begin{corollary}\label{cor2}
Let $K$ be a finite contractible flag $2$-complex which has no free edges.
Then $\dss(\Gamma_K)=3$ while $\gd(\Gamma_K)=2$. In these cases $\Gamma_K$ is finitely
presented (see \cite{BB}; an explicit presentation is given in \cite{DL}) and $FP_\infty$ (c.f. \cite{BB}).
\end{corollary}

The hypotheses on $K$ ensure that it cannot admit a piecewise Euclidean metric of
nonpositive curvature. The ``no free edges'' condition implies that any such metric
would have the geodesic extension property (see Proposition 5.10 of \cite{BH}), while
contractibility implies that the metric would be (globally) CAT(0). Together these 
condition force the metric to be unbounded, contradicting $K$ finite. 
Zeeman's dunce hat $D$ (a $2$-simplex $\Delta(a,b,c)$ with oriented edges $(a,b)$,
$(a,c)$ and $(b,c)$ all identified) is contractible with no free edges. Thus, for example, any flag
triangulation of $D$ will satisfy the hypotheses of Corollary \ref{cor2}.

\section{Refinement of Theorem \ref{main}}

\begin{figure}[ht]
\begin{center}
\includegraphics[width=13cm]{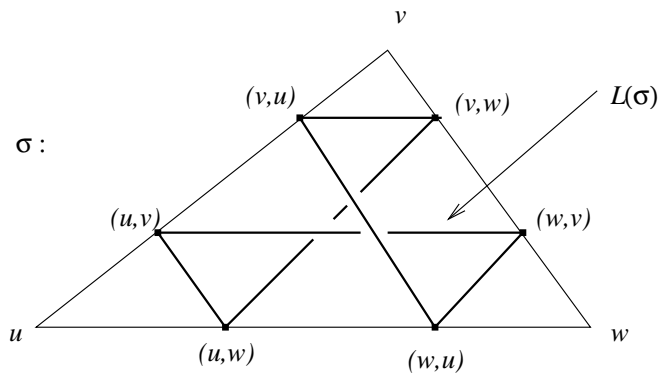}
\end{center}
\vglue -15pt
\caption{Defining $L(K)$ -- the contribution from a 2-simplex $\sigma$}\label{fig1}
\end{figure}

As in the previous section, let $K$ be a finite connected $2$-dimensional flag complex.
Suppose also that $K$ is given a piecewise Euclidean metric. We associate to $K$ a
metric graph $L(K)$ as follows. The vertex set of $L(K)$ is defined to be the set 
$\Cal E(K)$ of oriented edges of $K$. There is an edge of $L(K)$ between
$(u,u'),(v,v')\in\Cal E(K)$ if $(u,u')$ and $(v,v')$ are distinct edges of a common
$2$-simplex $\sig$ in $K$ and either $u=v$ or $u'=v'$. 
The length of such an edge of $L(K)$ is defined as the angle in $\sig$ between
the two sides in question.  The situation of a single simplex is illustrated in
Figure \ref{fig1}. Each simplex contributes to $L(K)$ a circle of length precisely
$2\pi$. These are identified pairwise along ``great 0-circles''
(pairs of antipodal points) according to the edge identifications between adjacent
$2$-simplexes in $K$. Note that $L(K)$ also contains, as a locally isometrically embedded
subgraph, the link of each vertex in $K$. Thus Theorem \ref{main} is a consequence of the following.

\begin{thm}\label{refined}
Let $K$ be a finite connected $2$-dimensional flag complex.
If $\Gamma_K$ acts properly by semi-simple isometries on a $2$-dimensional
(complete) CAT(0) space $X$, then $K$ admits a piecewise Euclidean metric such that 
$L(K)$ is a CAT(1) metric graph (or equivalently $L(K)$ 
contains no simple closed circuit of length strictly less than $2\pi$).
\end{thm}

\Proof We will use the fact, as pointed out in the remark following Lemma \ref{thetas},
that if $X$ is a complete CAT(0) metric space then its Tits boundary $\partial X$
is CAT(1). In particular, any locally geodesic circle in $\partial X$ has length
at least $2\pi$.

The metric on $K$ is defined exactly as in the proof of Theorem \ref{main}.
For any $2$-simplex $\sig$ in $K$, let $\Pi(\sig)$ denote the flat plane in $X$
associated with $\Gamma_\sig$. Then $L(\sig)$ is isometric to $\partial\Pi(\sig)$
by an isometry which associates a vertex $(u,v)$ of $L(\sig)$ with the boundary point 
$x_{(u,v)}^{+\infty}$. Since $\Pi(\sig)$ is a convex subspace of $X$ its boundary
embeds canonically as a locally convex subspace of $\partial X$. 
Thus we have a locally isometric map $L(\sig)\to\partial X$ for each
$2$-simplex $\sig$. (In fact, this map is isometric since $L(\sig)$ has 
diameter $\pi$). Combining these maps we have 
a map $\psi:L(K)\to\partial X$ which we claim is locally isometric. To check this, it
suffices to look at a neighbourhood of each vertex $(u,v)$ of $L(K)$. Such a
neighbourhood may be chosen to lie inside the subgraph $L(N(u,v))$,
where $N(u,v)$ denotes the union of all $2$-simplices containing the edge $(u,v)$. 
Now the $\Pi(\sig)$ for $\sig\in N(u,v)$ are mutually distinct planes all of which
lie in $\Min(x_{(u,v)})$. Since this minset is a convex subspace and, by the dimension
constraint on $X$, has the structure of $T\times\R$ where $T$ is an $\R$-tree,
it follows that $L(N(u,v))$ now embeds locally isometrically in $\partial X$.
This embedding is with respect to the intrinsic metric in the subgraph $L(N(u,v))$,
not the metric induced from $L(K)$ (these metrics do not agree if the subgraph is not
convex). However, locally these metrics do agree, showing that $\psi$ is a local isometry.

Since the local isometry $\psi$ maps simple closed circuits in $L(K)$ to locally
geodesic circles of the same length in $\partial X$ it now follows that all simple
closed circuits in $L(K)$ have length at least $2\pi$ and hence that $L(K)$ is CAT(1).

[Note that the hypothesis that $X$ be complete is actually artificial, since the same
result may be obtained with a little more effort by reworking the original proof of
Theorem \ref{main}. It is hypothesis (i) of Lemma \ref{thetas} which needs refining so as
to allow the case where $a_i=a_{i-1}a_{i+1}$, meaning that the three consecutive
elements lie in the same $\Z\times\Z$ subgroup but in such a way that one will always
have $\theta(a_{i-1},a_{i+1})=\theta(a_{i-1},a_i)+\theta(a_i,a_{i+1})$.]
\endproof

\begin{figure}[ht]
\begin{center}
\includegraphics[width=13cm]{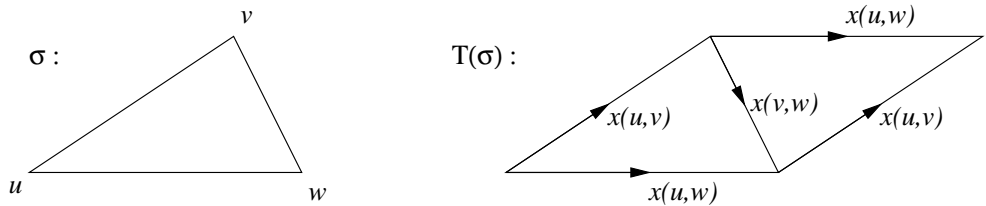}
\end{center}
\vglue -10pt
\caption{Building the complex $T(K)$ -- the torus $T(\sig)$ associated to a
$2$-simplex $\sig$}\label{fig2}
\end{figure}

Suppose now that $K$ is a finite \emph{simply-connected} piecewise Euclidean flag
complex. We now define a piecewise Euclidean complex $T(K)$ as follows. For each 
$2$-simplex $\sig$ in $K$ let $T(\sig)$ denote the the union of two isometric copies
of $\sig$ glued along their edges as shown in Figure \ref{fig2} so as to form a flat
torus. Oriented edges of $T(\sig)$ are labelled by elements of the generating set
$\Cal G(K)$ of $\Gamma_K$ as indicated in the figure. We now construct $T(K)$ by taking 
the union of all $T(\sig)$ for $\sig\in K$ and identifying edges (by isometries)
according to the labelling. Combinatorially, $T(K)$ is none other than the presentation
complex associated to the finite presentation for $\Gamma_K$ given in \cite{DL}
(Corollary 3) for the case $K$ simply-connected. Thus $\Gamma_K$ acts freely,
cocompactly and isometrically on the universal cover $\wtil{T(K)}$. Moreover, observe
that  the link of every vertex in $\wtil{T(K)}$, or rather the link of the unique
vertex in $T(K)$, is isometric to the metric graph $L(K)$. If $L(K)$ is CAT(1) then
$T(K)$ is non-positively curved and $\wtil{T(K)}$ is CAT(0). We therefore have the 
following Corollary to Theorem \ref{refined}.

\begin{corollary}
Let $K$ be a finite simply-connected $2$-dimensional flag complex.
Then the group $\Gamma_K$ acts properly by semi-simple isometries (or even freely
co-compactly) on a $2$-dimensional (complete) CAT(0) space $X$ if and only
if $K$ admits a piecewise Euclidean metric such that  $L(K)$ is a CAT(1) metric graph.
\end{corollary}

\begin{figure}[ht]
\begin{center}
\includegraphics[width=13cm]{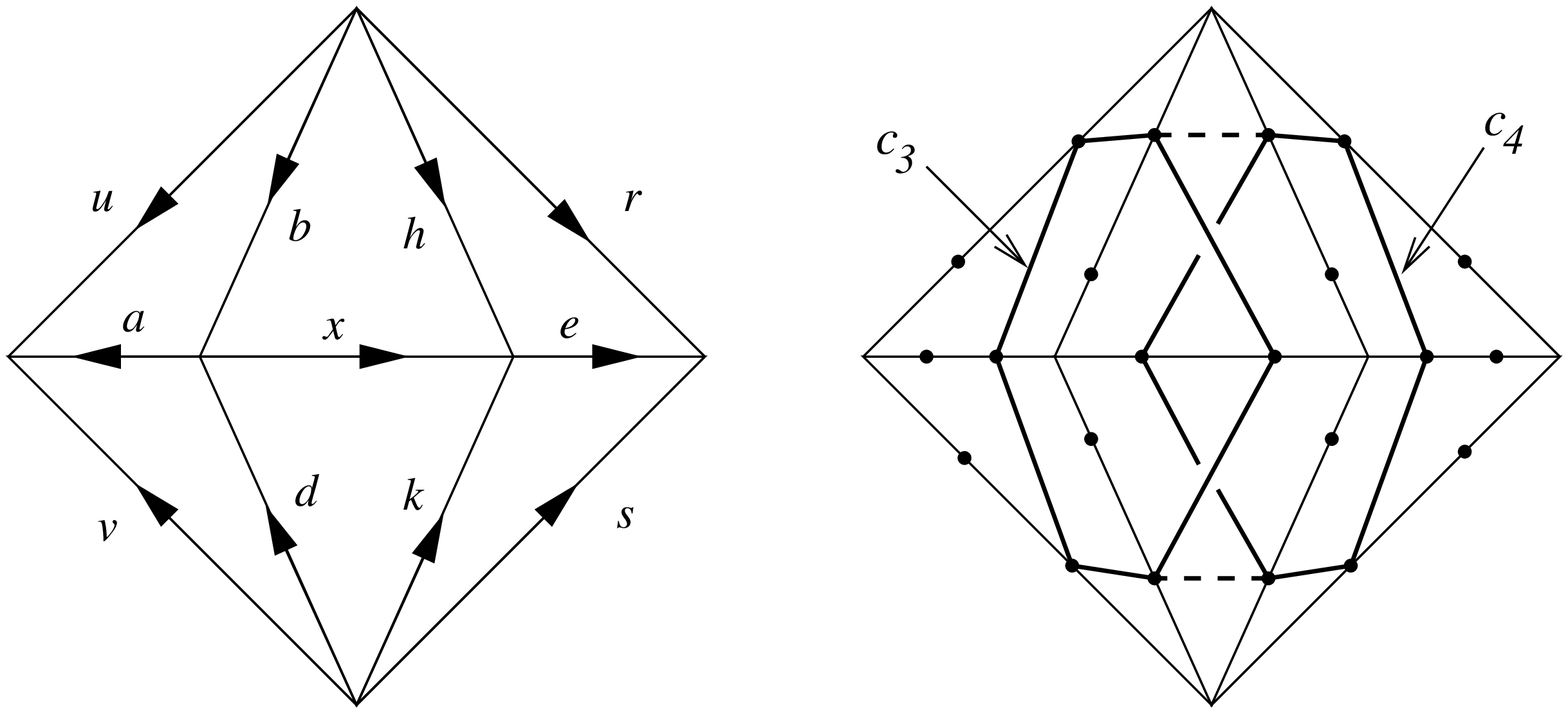}
\end{center}
\vglue -30pt
\caption{The flag complex $K_0$}\label{fig3}
\end{figure}

We finish with an example which shows that Theorem \ref{refined} really is
stronger than Theorem \ref{main}. The 2-dimensional flag complex $K_0$ shown in Figure
\ref{fig3} clearly admits a CAT(0) piecewise Euclidean metric, however we now show that
there is no piecewise Euclidean metric on $K_0$ such that $L_0=L(K_0)$ is a CAT(1)
graph. Suppose that $K_0$ did admit such a metric. In particular, $K_0$ would have to be
locally CAT(0), and hence CAT(0), with respect to this metric. We now consider the
lengths of the following circuits in $L_0$ (we describe a circuit by giving the
sequence of vertices of $L_0$ -- or rather, oriented edges of $K_0$ -- through which
it passes):   
\[
\begin{aligned}
c_1 &=  (b,h,r,e,s,k,d,v,a,u,b)\\
c_2 &= (b,x^{-1},d,k,x,h,b)\\
c_3 &= (b,x^{-1},d,v,a,u,b)\\
c_4 &= (h,r,e,s,k,x,h)
\end{aligned}
\]
The circuit $c_1$ has length exactly $2\pi$, since on the one hand it is a
circuit in the CAT(1) link $L_0$, so at least $2\pi$, while on the other hand its
length is exactly the angle sum of the quadrilateral in  $K_0$ with sides
$u,v,s,r$, so at most $2\pi$. The circuit $c_2$ has length $2\pi$ either by similar reasoning
or by simply noting that its length is the total angle sum of two Euclidean triangles.
But we also have
\[
\ell(c_3) +\ell(c_4) = \ell(c_1)+\ell(c_2) - 2d_{L_0}(b,h) - 2d_{L_0}(d,k) < 4\pi\,.
\]
Therefore at least one of the circuits $c_3$ or $c_4$ is strictly shorter than
$2\pi$, contradicting $L_0$ CAT(1).

By Theorem \ref{refined}, we have that $\dss(\Gamma_{K_0})=3$ while
$\gd(\Gamma_{K_0})=2$.

Let $L$ denote the subcomplex of $K_0$ which is the graph consisting of the three edges
labelled $a,x,e$ in Figure \ref{fig3}. Then $K_0$ is just the suspension of $L$. 
If $C(L)$ denotes the simplicial cone over $L$, then $\Gamma_{C(L)}$ is canonically isomorphic
to $A_L$. Thus $\Gamma_{K_0}$ is isomorphic to the double of $A_L$ over the normal subgroup
$\Gamma_L<A_L$:
\[
\Gamma_{K_0}\cong A_L\mathrel{\mathop{\star}\limits_{\Gamma_L}} A_L\,.
\]
Noel Brady \cite{Bra} has a proof that $\Gamma_{K_0}$ has a
Dehn function which is at least cubic. It follows (see \cite{BH}, for instance) 
that $\Gamma_{K_0}$ cannot act properly
cocompactly on any CAT(0) space, i.e: it is not a CAT(0) group. 
Without reproducing Brady's argument on the Dehn function we can still
deduce that this example is not a CAT(0) group as follows.
By the Exercises on page 499 of \cite{BH}, the group  $\Gamma_{K_0}$ may be re-expressed as 
the trivial HNN-extension of $A_L$ over $\Gamma_L$, and as such is a CAT(0) group
only if $\Gamma_L$ is quasi-isometrically embedded in $A_L$.
However, $\Gamma_L$ is sufficiently distorted in $A_L$ for this not to be the case.
Write $A_L=\<u_1,u_2,u_3,u_4 | [u_1,u_2]=[u_2,u_3]=[u_3,u_4]=1\>$, and observe that
$\Gamma_L\cong F_3$, freely generated by the elements $a=u_2^{-1}u_1$, $x=u_2^{-1}u_3$ and $e=u_3^{-1}u_4$.
For each $N\in\N$, the freely reduced word
$w_N:=(ax^{N}ex^{-N})^N$, of length $2(N^2+N)$, represents an element of the free group
$\Gamma_L$ which is equal in $A_L$ to
\[
\begin{aligned}
((u_2^{-1}u_1)u_2^{-N}u_3^N(u_3^{-1}u_4)u_3^{-N} u_2^N)^N &=(u_2^{-N}(u_2^{-1}u_1.u_3^{-1}u_4)u_2^N)^N\\
&= u_2^{-N}(u_2^{-1}u_1.u_3^{-1}u_4)^Nu_2^N\,.
\end{aligned}
\]
The latter word in the generators of $A_L$ has length $6N$, implying that $\Gamma_L$ is at least
quadratically distorted in $A_L$, and so is not quasi-isometrically embedded.

\Addresses\recd

\end{document}